


\magnification=\magstephalf
\frenchspacing

\parindent15pt

\abovedisplayskip4pt plus2pt
\belowdisplayskip4pt plus2pt
\abovedisplayshortskip2pt plus2pt
\belowdisplayshortskip2pt plus2pt

\font\twbf=cmbx10 at12truept
 at12truept
 at12truept

\font\sc=cmcsc10

\font\ninerm=cmr9 at9truept
\font\nineit=cmti9 at9truept
\font\ninesy=cmsy9 at9truept
\font\ninei=cmmi9 at9truept
\font\ninebf=cmbx9 at9truept

\font\sevenrm=cmr7 at7truept
 at7truept
\font\seveni=cmmi7 at7truept
\font\sevensy=cmsy7 at7truept

\font\fivenrm=cmr5 at5truept
\font\fiveni=cmmi5 at5truept
\font\fivensy=cmsy5 at5truept

\def\nine{%
\textfont0=\ninerm \scriptfont0=\sevenrm \scriptscriptfont0=\fivenrm
\textfont1=\ninei \scriptfont1=\seveni \scriptscriptfont1=\fiveni
\textfont2=\ninesy \scriptfont2=\sevensy \scriptscriptfont2=\fivensy
\textfont3=\tenex \scriptfont3=\tenex \scriptscriptfont3=\tenex
\def\rm{\fam0\ninerm}%
\textfont\itfam=\nineit
\def\it{\fam\itfam\nineit}%
\textfont\bffam=\ninebf
\def\bf{\fam\bffam\ninebf}%
\normalbaselineskip=11truept
\setbox\strutbox=\hbox{\vrule height8truept depth3truept width0pt}%
\normalbaselines\rm}

\hsize30truecc
\vsize44truecc
\nopagenumbers

\def\luz#1{\luzno#1?}
\def\luzno#1{\ifx#1?\let\next=\relax\yyy
\else \let\next=\luzno#1\xxx\fi\next}
\def\sp#1{\def\xxx{\kern1.7pt}\def\yyy{\kern-1.7pt}\luz{#1}}
\def\spa#1{\def\xxx{\kern1pt}\def\yyy{\kern-1pt}\luz{#1}}

\newcount\beg
\newbox\aabox
\newbox\atbox
\newbox\fpbox
\def\abbrevauthors#1{\setbox\aabox=\hbox{\sevenrm\uppercase{#1}}}
\def\abbrevtitle#1{\setbox\atbox=\hbox{\sevenrm\uppercase{#1}}}
\long\def\pag{\beg=\pageno
\def\leftheadline{\noindent\rlap{\nine\folio}\hfil\copy\aabox\hfil}
\def\rightheadline{\noindent\hfill\copy\atbox\hfill\llap{\nine\folio}}
\def\phead{\setbox\fpbox=\hbox{\sevenrm
************************************************}%
\noindent\vbox{\sevenrm\baselineskip9pt\hsize\wd\fpbox%
\centerline{***********************************************}

\centerline{BANACH CENTER PUBLICATIONS, VOLUME **}

\centerline{INSTITUTE OF MATHEMATICS}

\centerline{POLISH ACADEMY OF SCIENCES}

\centerline{WARSZAWA 19**}}\hfill}
\footline{\ifnum\beg=\pageno \hfill\nine[\folio]\hfill\fi}
\headline{\ifnum\beg=\pageno\phead
\else
\ifodd\pageno\rightheadline \else \leftheadline \fi
\fi}}

\newbox\tbox
\newbox\aubox
\newbox\adbox
\newbox\mathbox

\def\title#1{\setbox\tbox=\hbox{\let\\=\cr
\baselineskip14pt\vbox{\twbf\tabskip 0pt plus15cc
\halign to\hsize{\hfil\ignorespaces \uppercase{##}\hfil\cr#1\cr}}}}

\newbox\ab
\setbox\ab=\vbox{\vglue18pt}

\def\author#1{\setbox\aubox=\hbox{\let\\=\cr
\nine\baselineskip12pt\vbox{\tabskip 0pt plus15cc
\halign to\hsize{\hfil\ignorespaces \uppercase{\spa{##}}\hfil\cr#1\cr}}}%
\global\setbox\ab=\vbox{\unvbox\ab\box\aubox\vskip8pt}}

\def\address#1{\setbox\adbox=\hbox{\let\\=\cr
\nine\baselineskip12pt\vbox{\it\tabskip 0pt plus15cc
\halign to\hsize{\hfil\ignorespaces {##}\hfil\cr#1\cr}}}%
\global\setbox\ab=\vbox{\unvbox\ab\box\adbox\vskip16pt}}

\def\mathclass#1{\setbox\mathbox=\hbox{\footnote{}{1991 {\it Mathematics
Subject
Classification}\/: #1}}}

\long\def\maketitlebcp{\pag\unhbox\mathbox
\footnote{}{The paper is in final form and no version
of it will be published elsewhere.}
\vglue7truecc
\box\tbox
\box\ab
\vskip8pt}

\long\def\abstract#1{{\nine{\bf Abstract.}
#1

}}

\def\section#1{\vskip-\lastskip\vskip12pt plus2pt minus2pt
{\bf #1}}

\long\def\th#1#2#3{\vskip-\lastskip\vskip4pt plus2pt
{\sc #1} #2 {\it #3}\vskip-\lastskip\vskip4pt plus2pt}

\long\def\remar#1{\vskip-\lastskip\vskip4pt plus2pt
\sp{#1}\ \ignorespaces}

\def\Proof{\vskip-\lastskip\vskip4pt plus2pt
\sp{Proo{f.}\ }\ignorespaces}

\def\endproof{\nobreak\kern5pt\nobreak\vrule height4pt width4pt depth0pt
\vskip4pt plus2pt}

\newbox\refbox
\newdimen\refwidth
\long\def\ref#1#2{{\nine
\setbox\refbox=\hbox{\nine[#1]}\refwidth\wd\refbox\advance\refwidth by 12pt%
\def\textindent##1{\indent\llap{##1\hskip12pt}\ignorespaces}
\vskip24pt plus4pt minus4pt
\centerline{\bf References}
\vskip12pt plus2pt minus2pt
\parindent=\refwidth
#2

}}

\def\footnoterule{\kern -3pt \hrule width 4truecc \kern 2.6pt}

\catcode`@=11
\def\vfootnote#1%
{\insert\footins\bgroup\nine\interlinepenalty\interfootnotelinepenalty%
\splittopskip\ht\strutbox\splitmaxdepth\dp\strutbox\floatingpenalty\@MM%
\leftskip\z@skip\rightskip\z@skip\spaceskip\z@skip\xspaceskip\z@skip%
\textindent{#1}\footstrut\futurelet\next\fo@t}
\catcode`@=12

\mathclass{Primary 14C15, 14H10; Secondary 14H42.}

\abbrevauthors{C.F.~Faber}
\abbrevtitle{Intersection-theoretical computations}

\title{Intersection-theoretical computations on
${{\overline {\cal M}}_{\lowercase{g}}}$}

\author{Carel\ Faber}
\address{Faculteit der Wiskunde en Informatica\\
Universiteit van Amsterdam\\
Plantage Muidergracht 24\\
1018 TV \ Amsterdam\\
Netherlands\\
E-mail: faber@fwi.uva.nl}

\maketitlebcp

\input amssym.def
\input amssym.tex
\font\tenfrak=eufm10
\newfam\frakfam
\textfont\frakfam=\tenfrak
\def\frak{\fam\frakfam\tenfrak}

\def\no{}

\def\pn{\par}

\def\Z   {{\Bbb Z}}

\def\P   {{\Bbb P}}
\def\Q   {{\Bbb Q}}
\def\C   {{\Bbb C}}

\def\mbar#1{{\overline{\cal M}}_{#1}}

\def\mgbar{{\overline {\cal M}}_g}

\def\la#1{\lambda_{#1}}
\def\deg{{\rm deg\ }}
\def\ll{\lambda}
\def\d#1{\delta_{#1}}
\def\D#1{\Delta_{#1}}
\def\ka#1{\kappa_{#1}}

\section{Introduction.}
In this paper we explore several concrete problems, all more or less
related to the intersection theory of the moduli space of (stable) curves,
introduced by Mumford [Mu 1].
\pn In \S1 we only intersect divisors with curves. We find a collection of
necessary conditions for ample divisors, but the question whether these
conditions are also sufficient is very much open.
\pn The other sections are concerned with moduli spaces of curves of low genus,
but we use the ring structure of the Chow ring. In \S\S2, 3 we find necessary
conditions for very ample divisors on $\mbar2$ and $\mbar3$.
\pn The intersection numbers of the kappa-classes are the subject of the Witten
conjecture, proven by Kontsevich. In \S4 we show how to compute these numbers
for $g=3$ within the framework of algebraic geometry.
\pn Finally, in \S5 we compute $\ll^9$ on $\mbar4$. This also gives the
value of $\la{g-1}^3$ (for $g=4$), which is relevant for counting
curves of higher genus on manifolds [BCOV]. Another corollary is a different
computation of the class of the Jacobian locus in the moduli space of
4-dimensional principally polarized abelian varieties; in a sense this gives
also a different proof that the Schottky locus is irreducible in
dimension 4.
\pn {\sl Acknowledgement.\/} I would like to thank Gerard van der Geer for very
useful discussions in connection with \S5. This research has been made possible
by a fellowship of the Royal Netherlands Academy of Arts and Sciences.

\section{1. Necessary conditions for ample divisors on $\mgbar$ .}
Let $g\ge2$ be an integer and put $h=[g/2]$. Cornalba and Harris [C-H]
determined which divisors on $\mgbar$ of the form $a\lambda-b\delta$
are ample: this is the case if and only if $a>11b>0$. Divisors of this
form are numerically effective (nef) if and only if $a\ge11b\ge0$.
(More generally, the ample cone is the interior of the nef cone and the
nef cone is the closure of the ample cone ([Ha], p.~42)).
Here $\delta=\sum_{i=0}^{h}\delta_i$ with $\delta_i=[\Delta_i]$
for $i\neq1$ and $\delta_1={1\over2}[\Delta_1]$.
\pn Arbarello and Cornalba [A-C] proved that the $h+2$ divisors
$\lambda,\delta_0,\delta_1,\dots,\delta_h$ form for $g\ge3$ a $\Z$-basis of
$Pic(\mgbar)$
(the Picard group of the moduli functor),
using the results of Harer and
Mumford (we work over $\C$).
As pointed out in [C-H] it would be interesting to
determine the nef cone in $Pic(\mgbar)$ for $g\ge3$. (For $g=2$
the answer is given by the result of [C-H], because of the relation
$10\lambda-\delta_0-2\delta_1=0$.)
\pn In [Fa 1], Theorem 3.4, the author determined the nef cone for $g=3$.
The answer is: $a\lambda-b_0\delta_0-b_1\delta_1$ is nef on $\mbar 3 $
if and only if $2b_0\ge b_1\ge0$ and $a-12b_0+b_1\ge0$.
That a nef divisor necessarily satisfies these inequalities, follows from
the existence of one-dimensional families of curves for which
$(\deg\lambda,\deg\delta_0,\deg\delta_1)$ equals $(1,12,-1)$
resp.~$(0,-2,1)$ resp.~$(0,0,-1)$.
Such families are easily constructed: for the first family, take a
simple elliptic pencil and attach it to a fixed one-pointed curve of genus
$2$;
for the second family, take a $4$-pointed rational
curve with one point moving and attach a fixed two-pointed curve of genus
$1$ to two of the points and identify the two other points;
for the third family, take a $4$-pointed rational
curve with one point moving and attach two fixed one-pointed curves of
genus 1 to two of the points and identify the two other points.
\pn That a divisor on $\mbar 3$
satisfying the inequalities is nef, follows once we show
that $\lambda$, $12\lambda-\delta_0$ and $10\lambda-\delta_0-2\delta_1$
are nef. It is well-known that $\lambda$ is nef. Using induction on the genus
one shows that $12\lambda-\delta_0$ is nef:
on $\mbar{1,1}$ it vanishes; for $g\ge2$, writing
$12\lambda-\delta_0=\kappa_1+\sum_{i=1}^h\delta_i$ one sees that
$12\lambda-\delta_0$ is positive on every one-dimensional family of curves
where the generic fiber has at most nodes of type $\delta_0$; if on the other
hand the generic fiber has a node of type $\delta_i$ for some $i>0$, one
partially normalizes the family along a section of such nodes and uses the
induction hypothesis (cf.~the proof of Proposition 3.3 in [Fa 1], which
unfortunately proves the result only for $g=3$).
Finally, the proof that $10\lambda-\delta_0-2\delta_1$ is nef on $\mbar 3$
is ad hoc (see the proof of Theorem 3.4 in [Fa 1]).
\pn All we do in this section is come up with a couple of one-dimensional
families of stable curves for which we compute the degrees of the basic
divisors. The naive hope is that at least some of these families are
extremal (cf.~[C-H], p.~475), but the author hastens to add that there
is at present very little evidence to support this.
\pn The method of producing families is a
very simple one: we start out trying to
write down all the families for which the generic fiber has $3g-4$ nodes.
This turns out to be a bit complicated. However, the situation greatly
simplifies as soon as one realizes that the only one-dimensional
moduli spaces of stable pointed curves are $\mbar{0,4}$ and
$\mbar{1,1}$: for the computation of the basic divisor classes on these
families, one only needs to know the genera of the pointed curves attached
to the moving $4$-pointed rational curve resp.~the moving one-pointed curve
of genus $1$ as well as the types of the nodes one gets in this way. In other
words, the fixed parts of the families can be taken to be general.
\pn We now consider the various types of families obtained in this way
and compute on each family the degrees of the basic divisor classes.
Each family gives a necessary condition for the divisor $a\lambda
-\sum_{i=0}^hb_i\delta_i$ to be nef. In order to write this condition,
it will be convenient to define $\delta_i=\delta_{g-i}$ and $b_i=b_{g-i}$
for $h<i<g$.
\pn A) In the case of $\mbar{1,1}$, there is very little choice: we can only
attach a (general) one-pointed curve of genus $g-1$. Taking a simple
elliptic pencil for the moving part, we get---as is well-known---the
following degrees: $\deg\lambda=1$, $\deg\delta_0=12$, $\deg\delta_1=-1$
and $\deg\delta_i=0$ for $1<i\le h$. This gives the necessary condition
$a-12b_0+b_1\ge0$.
\pn B) The other families are all constructed from a $4$-pointed smooth
rational curve with one of the points moving and the other three fixed;
when the moving point meets one of the fixed points, the curve breaks
up into two 3-pointed smooth rational curves glued at one point. We have
to examine the various ways of attaching general curves to this
$4$-pointed rational curve. E.g., one can attach one curve, necessarily
$4$-pointed and of genus $g-3$. All nodes are of type $\delta_0$ and
the 3 degenerations have an extra such node. Therefore
$\deg\delta_0=-4+3=-1$, while the other degrees are zero; one obtains the
necessary condition $b_0\ge0$.
\pn C) Now attach a $3$-pointed curve of genus $i$ and a 1-pointed
curve of genus $j\ge1$, with $i+j=g-2$. One checks $\deg\delta_0
=-3+3=0$ and $\deg\delta_j=-1$, the other degrees vanish. One obtains
$b_j\ge0$ for $j\ge1$.
Thus all $b_i$ are non-negative for a nef divisor.
\remar{Remark.} If one uses the families above, one simplifies the proof
of Theorem 1 in [A-C] a little bit.
\vskip4pt plus2pt

\pn D) If we attach two-pointed curves of genus $i\ge1$ and $j\ge1$, with
$i+j=g-2$, we find $\deg\delta_0=-4+2=-2$ and $\deg\delta_{i+1}=1$. So for
$2\le k\le h$ we find the condition $2b_0-b_k\ge0$.
\pn E) Attaching a two-pointed curve of genus $i$ and two one-pointed curves
of genus $j$ and $k$, with $i,j,k\ge1$ and $i+j+k=g-1$, we find that two
of the degenerations have an extra node of type $\delta_0$ while the third
has an extra node of type $\delta_{j+k}$. Therefore $\deg\delta_0=-2+2=0$.
It is cumbersome to distinguish the various cases that occur for the
other degrees, but is also unnecessary: one may simply write the
resulting necessary condition in the form $b_j+b_k-b_{j+k}\ge0$, for
$j,k$ with $1\le j\le k$ and $j+k\le g-2$.
\pn F) Attaching 4 one-pointed curves of genera $i,j,k,l\ge1$, with
$i+j+k+l=g$, we get the necessary condition
$b_i+b_j+b_k+b_l-b_{i+j}-b_{i+k}-b_{i+l}\ge0$.
\pn G) If we identify two of the 4 points to each other and attach a
two-pointed curve of genus $g-2$ to the remaining two points, we obtain
the necessary condition $2b_0-b_1\ge0$.
\pn H) As in G), but now we attach 1-pointed curves of genera $i,j\ge1$
to the remaining two points, with $i+j=g-1$. The resulting condition
is $b_i+b_j-b_1\ge0$.
\pn The only other possibility is to identify the first with the second
and the third with the fourth point. This gives a curve of genus 2, so
this is irrelevant.
We have proven the following theorem.
\th{Theorem}{1.}{Assume $g\ge3$. A numerically effective divisor
$a\lambda-\sum_{i=0}^hb_i\delta_i$ in $Pic(\mbar g)$ satisfies
the following conditions:
\item{a)} $a-12b_0+b_1\ge0\quad;$
\item{b)} for all $j\ge1$,
$$2b_0\ge b_j\ge 0\quad;$$
\item{c)} for all $j,k$ with $1\le j\le k$ and $j+k\le g-1$,
$$b_j+b_k\ge b_{j+k}\quad;$$
\item{d)} for all $i,j,k,l$ with $1\le i\le j\le k\le l$ and
$i+j+k+l=g$,
$$b_i+b_j+b_k+b_l\ge b_{i+j}+b_{i+k}+b_{i+l}\quad.$$}

\no Here $b_i=b_{g-i}$ for $h<i<g$, as before.
The conditions in the theorem are somewhat redundant. E.g., it is
easy to see that condition (c) implies the non-negativity of
the $b_i$ with $i\ge1$.
\pn As we have seen, the conditions in the theorem are sufficient for
$g=3$. The proof proceeded by determining the extremal rays of the cone
defined by the inequalities and analyzing the (three) extremal rays separately.
It may therefore be of some interest to find (generators for) the extremal
rays of the cone in the theorem. We have done this for low genus:
$$
\eqalign{&g=4:\qquad \cases{\lambda\cr
12\lambda-\delta_0\cr
10\lambda-\delta_0-2\delta_1\cr
10\lambda-\delta_0-2\delta_1-2\delta_2\cr
21\lambda-2\delta_0-3\delta_1-4\delta_2\cr}
\cr
&g=5:\qquad\cases{\lambda\cr
12\lambda-\delta_0\cr
10\lambda-\delta_0-2\delta_1-\delta_2\cr
10\lambda-\delta_0-2\delta_1-2\delta_2\cr
32\lambda-3\delta_0-4\delta_1-6\delta_2\cr}
\cr
&g=6:\qquad\cases{\lambda\cr 12\lambda-\delta_0\cr
10\lambda-\delta_0-2\delta_1-2\delta_2\cr
10\lambda-\delta_0-2\delta_1-2\delta_3\cr
10\lambda-\delta_0-2\delta_1-2\delta_2-2\delta_3\cr
32\lambda-3\delta_0-4\delta_1-6\delta_2-6\delta_3\cr
98\lambda-9\delta_0-10\delta_1-16\delta_2-18\delta_3\cr}
\cr}$$

\pn Unfortunately, we have not been able to discover a general pattern.
(There are 10 extremal divisors for $g=7$,
20 extremal divisors for $g=8$ and
21 extremal divisors for $g=9$.)
It is easy to see that $\lambda$, $12\lambda-\delta_0$ and
$10\lambda-2\delta+\delta_0$ are extremal in every genus.
It should be interesting to know the answer to the following question.
\th{Question.}{}{\item{a)} Is $10\lambda-2\delta+\delta_0$ nef for all $g\ge4$?
\item{b)} Are the conditions in the theorem sufficient??}

\no Note that an affirmative answer to the first question implies the
result of [C-H] mentioned above, since $12\lambda-\delta_0$ is nef.
Note also that a divisor satisfying the conditions in the theorem is
non-negative on every one-dimensional family of curves
whose general member is smooth. This follows easily from
[C-H, (4.4) and Prop.~(4.7)].
(I would like to thank Maurizio Cornalba for reminding me
of these results.)

\section{2. Necessary conditions for very ample divisors on $\mbar2$ .}
We know which divisors on $\mbar2$ are ample: it is easy to see that
$\lambda$ and $\delta_1$ form a $\Z$-basis of the functorial Picard group
$Pic(\mbar2)$; then $a\lambda+b\delta_1$ is ample if and only if
$a>b>0$, as follows from the relation $10\lambda=\delta_0+2\delta_1$ and
the fact that $\lambda$ and $12\lambda-\delta_0$ are nef.
\pn Therefore it might be worthwhile to study which divisors are very ample
on the {\it space\/} $\mbar2$.
Suppose that $D=a\lambda+b\delta_1$ is a very ample divisor. Then for every
$k$-dimensional subvariety $V$ of $\mbar2$ the intersection product
$D^k\cdot[V]$ is a positive integer, the degree of $[V]$ in the embedding
of $\mbar2$ determined by $|D|$. We work this out for the subvarieties
that we know; we use Mumford's computation [Mu~1] of the Chow ring (with
$\Q$-coefficients) of $\mbar 2$.
The result may be formulated as follows:
$$
A^*(\mbar2)=\Q[\lambda,\delta_1]/(\lambda(\lambda+\delta_1),\lambda^2
(5\lambda-\delta_1)).
$$
The other piece of information we need is on p.~324 of [Mu 1]:
$\lambda^3={1\over2880}p$. However, one should realize that the identity
element in $A^*(\mbar2)$ is $[\mbar2]_Q={1\over2}[\mbar2]$, which means
that
$$\lambda^3\cdot[\mbar2]={1\over1440}\quad.$$
Therefore
$$D^3\cdot[\mbar2]={{a^3+15a^2b-15ab^2+5b^3}\over1440}\quad.$$
One of the requirements is therefore that the integers $a$ and $b$ are
such that the expression above is an integer. It is not hard to see
that this is the case if and only if
$$60|a\qquad\hbox{and}\qquad12|b\quad.$$
It turns out that these conditions imply that $D^2\cdot[\Delta_0]$
and $D^2\cdot[\Delta_1]$ are integers. Also $D^2\cdot4\lambda$ is an integer,
but $D^2\cdot2\lambda$ is an integer if and only if $8|(a+b)$. Therefore,
if for some integer $k$ the class $(4k+2)\lambda$ is the fundamental class
of an effective $2$-cycle, then a very ample $D$ satisfies $8|(a+b)$.
We don't know whether such a $k$ exists; clearly, $20\lambda=
[\Delta_0]+[\Delta_1]$ is effective; the fundamental class of the bi-elliptic
divisor turns out to be $60\lambda+3\Delta_1$.
\pn Turning next to one-dimensional subvarieties, the conditions $60|a$ and
$12|b$ imply that $D\cdot[\Delta_{00}]$ and $D\cdot[\Delta_{01}]$
are integers as well.
\th{Proposition}{2.}{A very ample divisor $a\lambda+b\delta_1$ on the
moduli space $\mbar2$ satisfies the following conditions:
\item{a)} $a,b\in\Z$ and $a>b>0\quad;$
\item{b)} $60|a$ and $12|b\quad.$}

\th{Corollary}{3.}{The degree of a
projective embedding of $\mbar2$ is at least $516$.}

\Proof We need to determine for which $a$ and $b$ satisfying
the conditions in the proposition the expression $5(b-a)^3+6a^3$ attains
its minimum value. Clearly this happens exactly for $b=12$ and $a=60$.
If $60\lambda+12\delta_1$ is very ample, the degree of $\mbar2$ in
the corresponding embedding is $(5(b-a)^3+6a^3)/1440=516$.
\vskip4pt plus2pt

\remar{Remark.} It is interesting to compare
the obtained necessary conditions
with the explicit descriptions of $\mbar2$ given by Qing Liu ([Liu]).
The computations we have done
(in characteristic $0$) indicate that $60\lambda+60\delta_1$ maps
$\mbar2$ to a copy of $\frak X$ (loc.~cit., Th\'eor\`eme 2),
that $60\lambda+36\delta_1$ maps $\mbar2$ to the blowing-up of $\frak X$
with center ${\cal J}_{\Q}$ (loc.~cit., Corollaire 3.1)
and that $60\lambda+48\delta_1$ is very ample, realizing $\mbar2$ as the
blowing-up of $\frak X$ with center the ideal generated by
$I_4^3$, $J_{10}$, $H_6^2$ and $I_4^2H_6$ (loc.~cit., Corollaire 3.2).
\vskip4pt plus2pt

\section{3. Necessary conditions for very ample divisors on $\mbar3$ .}
In this section we compute necessary conditions for very ample
divisors on the moduli space $\mbar3$. As we mentioned in \S1,
a divisor $D=a\ll-b\d0-c\d1\in Pic(\mbar3)$ with $a,b,c\in\Z$ is
ample if and only if $a-12b+c>0$ and $2b>c>0$. The necessary conditions
for very ample $D$
are obtained as in \S2: for a $k$-dimensional subvariety $V$ of $\mbar3$,
the intersection product $D^k\cdot[V]$ should be an integer. We use the
computation of the Chow ring of $\mbar3$ in [Fa 1]. The computations are
more involved than in the case of genus 2; also, we know the
fundamental classes of more subvarieties.
\pn First we look at the degree of $\mbar3$:
$$
\eqalign{
D^6=(a\ll-b\d0-c\d1)^6
&=\textstyle{1\over90720}a^6-{1\over576}a^4c^2-{1\over18}a^3b^3
+{1\over48}a^3bc^2+{35\over3456}a^3c^3\cr
&\qquad\textstyle
+{5\over8}a^2b^2c^2-{43\over96}a^2bc^3+{13\over512}a^2c^4
+{203\over20}ab^5-{145\over12}ab^3c^2\cr
&\qquad\textstyle
+{25\over4}ab^2c^3-{31\over48}abc^4+{149\over7680}ac^5
-{4103\over72}b^6+{55}b^4c^2
\cr
&\qquad\textstyle
-{505\over18}b^3c^3+{65\over16}b^2c^4-{91\over384}bc^5+{5\over1024}c^6
\ ,
\cr}
$$
as follows from [Fa 1], p.~418.
The requirement that this is in $\Z_2$ implies, firstly, that $2|c$,
secondly, that $2|a$ {\sl and\/} $4|c$, thirdly, that $2|b$. Looking in
$\Q_3$ we get, firstly, that $3|a$, secondly, that $3|b$.
Modulo 5 we get $5|a$ or $5|(a+3b+c)$. Finally, working modulo 7
we find that $7|a$ should hold.
\pn Writing $a=42a_1$, $b=6b_1$ and $c=4c_1$, with $a_1,b_1,c_1\in\Z$, the
condition $D^6\cdot[\mbar3]\in\Z$ becomes $5|a_1$ or $5|(3a_1+2b_1+c_1)$.
Interestingly, unlike the case of genus 2, these conditions are not the
only necessary conditions we find.
\pn For instance, the condition $D^5\cdot\d0\in\Z$ translates in
$3|c_1$; then $[\Delta_1]=2\d1$ gives no further conditions; but the
hyperelliptic locus, with fundamental class $[{\cal H}_3]=
18\ll-2\d0-6\d1$, improves the situation modulo 5: necessarily
$5|(3a_1+2b_1+c_1)$. It follows that $D^5\cdot\ll$ is an integer,
so all divisors have integer-valued degrees.
\pn In codimension 2, writing $c_1=3c_2$ with $c_2\in\Z$, the condition
$D^4\cdot[\Delta_{01a}]\in\Z$ translates in
$$
5|a_1\qquad\hbox{or}\qquad5|c_2\qquad\hbox{or}\qquad5|(a_1+c_2)\qquad\hbox{or}
\qquad5|(a_1+3c_2)\ .
$$
The (boundary) classes
$[\D{00}]$, $[\D{01b}]$, $[\D{11}]$, $[\Xi_0]$, $[\Xi_1]$ and $[{\rm H}_1]$
([Fa 1], pp.~340 sqq.) give no further conditions.
\pn In codimension 3, the class $[(i)]=8[(i)]_Q$ forces $2|a_1$. Write
$a_1=2a_2$ with $a_2\in\Z$.
Somewhat surprisingly, the class $[{\rm H}_{01a}]=4\eta_0$ (loc.~cit., pp.~386,
388) gives the condition $5|(a_2+2c_2)$. Consequently, combining the
various conditions modulo 5, we obtain
$$
5|a_2\qquad\hbox{and}\qquad5|b_1\qquad\hbox{and}\qquad5|c_2\ .
$$
Finally, we checked that the 12 cycles in codimension 4 and the
8 cycles in codimension 5 (loc.~cit., pp.~346 sq.)
don't give extra conditions.
\th{Proposition}{4.}{A very ample divisor $a\ll-b\d0-c\d1$ on the
moduli space $\mbar3$ satisfies the following conditions:
\item{a)} $a,b,c\in\Z$ with $a-12b+c>0$ and $2b>c>0\quad;$
\item{b)} $420|a$ and $30|b$ and $60|c\quad.$}

\th{Corollary}{5.}
{The degree of a projective embedding of $\mbar3$ is at least
$$
650924662500=2^2\cdot3^2\cdot5^5\cdot7\cdot826571.
$$}

\Proof We need to minimize the expression given for the
degree of $\mbar3$ while fulfilling the conditions in the proposition.
Write $a=420A$, $b=30B$ and $c=60C$. One shows that in the cone given by
$7A-6B+C\ge0$ and $B\ge C\ge0$ the degree is minimal along the (extremal)
ray $(A,B,C)=(5x,7x,7x)$ (corresponding to $10\ll-\d0-2\d1$).
Comparing the value for $(A,B,C)=(5,7,7)$ with that for
$(A,B,C)=(2,2,1)$, one concludes $A\le5$, $B\le7$ and $C\le7$. This leaves
only a few triples in the interior of the cone; the minimum degree
is obtained for $(A,B,C)=(2,2,1)$, corresponding to $840\ll-60\d{}$.
\vskip4pt plus2pt

\remar{Remark.} In [Fa 1], Questions 5.3 and 5.4, we asked whether
the classes $X$ (resp.~$Y$) are multiples of classes of complete
subvarieties of $\mbar3$ of dimension 4 (resp.~3) having empty
intersection with $\D1$ (resp.~$\D0$). We still don't know the answers,
but we verified that $X$ and $-Y=504\la3$ are effective:
$$
\eqalign{
X&=\textstyle{1\over15}\d{00}+{1\over6}\d{01a}+{11\over15}\d{01b}
+8\d{11}+{3\over14}\xi_0+{48\over35}\xi_1+{40\over21}\eta_1\quad;\cr
-Y&=\textstyle{1\over2}[(a)]_Q+[(b)]_Q+[(c)]_Q+{11\over30}[(d)]_Q
+{2\over5}[(f)]_Q+2[(g)]_Q+{2\over3}\eta_0\quad.\cr
}
$$
(For the notations, see [Fa 1], pp.~343, 386, 388.)
\vskip4pt plus2pt

\section{4. Algebro-geometric calculation of the intersection numbers
of the tautological classes on $\mbar3$ .}
Here we show how to compute the intersection numbers of the classes
$\ka{i}$ $(1\le i\le6)$ on $\mbar3$ in an algebro-geometric setting.
These calculations were done originally in May 1990 to check the
genus 3 case of Witten's conjecture [Wi], now proven by Kontsevich [Ko].
We believe that there is still interest, though, in
finding methods within algebraic geometry that allow to compute the
intersection numbers of the kappa- or tau-classes. For instance, the
identity
$$K^{3g-2}=\langle\tau_{3g-2}\rangle=\langle\ka{3g-3}\rangle
={1\over(24)^g\cdot g!}
$$
(in cohomology) should be understood ([Wi], between (2.26) and (2.27)).

\pn In [Fa 1] the 4 intersection numbers of $\ka1$ and $\ka2$ were
computed; using the identity $\ka1=12\ll-\d0-\d1$, we can read these off
from Table 10 on p.~418:
$$
\textstyle
\ka1^6={176557\over107520}\quad,\quad\ka1^4\ka2={75899\over322560}\quad,
\quad\ka1^2\ka2^2={32941\over967680}\quad,\quad\ka2^3={14507\over2903040}\quad.
$$
To compute the other intersection numbers, we need to express the other
kappa-classes in terms of the bases introduced in [Fa 1]. The set-up is
as in [Mu 1], \S8 (and \S6): if $C$ is a stable curve of genus 3,
$\omega_C$ is generated by its global sections, unless
\item{a)} $C$ has 1 or 2 nodes of type $\d1$, in which case the
global sections generate the subsheaf of $\omega_C$ vanishing in these
nodes;
\item{b)} $C$ has 3 nodes of type $\d1$, i.e., $C$ is a $\P^1$
with 3 (possibly singular) elliptic tails, in which case
$\Gamma(\omega_C)$ generates the subsheaf of $\omega_C$ of sections
vanishing on the $\P^1$.

\noindent (See [Mu 1], p.~308.) Let $Z\subset{\overline {\cal C}}_3$ be the
closure of the locus of pointed curves with 3 nodes of type $\d1$ and with the
point lying on the $\P^1$. Working over ${\overline {\cal C}}_3-Z$ we get
$$
0\to{\cal F}\to\pi^*\pi_*
\omega_{{\overline {\cal C}}_3/\mbar3}\to I_{\D1^*}\cdot
\omega_{{\overline {\cal C}}_3/\mbar3}\to0
$$
with ${\cal F}$ locally free of rank 2. Working this out as in [Fa 1], p.~367
we get
$$
\eqalign{
0=c_3({\cal F})&=\pi^*\la3-K\cdot\pi^*\la2+K^2\cdot\pi^*\la1-K^3\cr
&\qquad-(\pi^*\la1-K)\cdot[\D1^*]_Q+i_{1,*}(K_1+K_2)
\cr}
$$
modulo $[Z]$. Multiplying this with $K$ and using that $\omega^2$ is trivial
on $[\D1^*]$, we get
$$0=K\cdot c_3({\cal F})=K\cdot\pi^*\la3-K^2\cdot\pi^*\la2
+K^3\cdot\pi^*\la1-K^4+*K\cdot[Z]\quad.\leqno(1)
$$
It is easy to see that $K^2\cdot[Z]=0$, so we also get
$$
\leqalignno{
0&=K^2\cdot\pi^*\la3-K^3\cdot\pi^*\la2
+K^4\cdot\pi^*\la1-K^5\quad,&(2)\cr
0&=K^3\cdot\pi^*\la3-K^4\cdot\pi^*\la2
+K^5\cdot\pi^*\la1-K^6\quad,&(3)\cr
0&=K^4\cdot\pi^*\la3-K^5\cdot\pi^*\la2
+K^6\cdot\pi^*\la1-K^7\quad.&(4)\cr
}
$$
Pushing-down to $\mbar3$ we get
$$
\leqalignno{
0&=4\la3-\ka1\la2+\ka2\la1-\ka3+N\cdot[(i)]_Q\quad,&(1')\cr
0&=\ka1\la3-\ka2\la2+\ka3\la1-\ka4\quad,&(2')\cr
0&=\ka2\la3-\ka3\la2+\ka4\la1-\ka5\quad,&(3')\cr
0&=\ka3\la3-\ka4\la2+\ka5\la1-\ka6\quad.&(4')\cr
}
$$
To get $\ka3$ from $(1')$ we use two things. Firstly, one computes
$Y=-504\la3$, as mentioned at the end of \S3. This follows since both
$Y$ and $\la3$ are in the one-dimensional subspace of $A^3(\mbar3)$
of classes vanishing on all subvarieties of $\D0$. The factor $-504$
is computed using $\ll^4=8\ll\la3$ or $\la3\cdot[(i)]_Q={1\over6}\ll^3
\cdot[(i)]_Q$.
Secondly, to compute $N$, one uses that $\ka3$ vanishes on the classes
$[(b)]_Q$, $[(c)]_Q$, $[(f)]_Q$, $[(g)]_Q$, $[(h)]_Q$ and $[(i)]_Q$.
This gives 6 relations in $N$ of which 3 are identically zero; the other
3 all imply $N=1$.
\pn The formulas above allow one to express the kappa-classes in terms
of the bases of the Chow groups given in [Fa 1].
We give the formula for $\ka3$ (from which the other formulas follow):
$$\eqalign{\textstyle
\ka3&\textstyle={1\over280}[(a)]_Q+{31\over840}[(b)]_Q+{19\over420}[(c)]_Q
+{1\over1260}[(d)]_Q+{1\over35}[(e)]_Q\cr
&\textstyle\qquad+{19\over840}[(f)]_Q+{29\over84}[(g)]_Q
+{11\over35}[(h)]_Q+{93\over35}[(i)]_Q+{11\over252}\eta_0\quad.
\cr
}
$$
This gives the following intersection numbers:
$$
\displaylines{
\textstyle
\ka1^3\ka3={4073\over161280}\quad,\quad\ka1\ka2\ka3={149\over40320}\quad,\quad
\ka3^2={131\over322560}\quad,\quad\ka1^2\ka4={2173\over967680}\quad,\cr
\textstyle
\ka2\ka4={971\over2903040}\quad,\quad\ka1\ka5={1\over5760}\quad,\quad
\ka6={1\over82944}
\quad.\cr
}
$$

\section{5. A few intersection numbers in genus 4.}
Kontsevich's proof of Witten's conjecture enables one to compute the
intersection numbers of the kappa-classes on the moduli space of
stable curves of arbitrary genus. There are many more intersection
numbers that one would like to know, see e.g.~[BCOV], (5.54) and end of
Appendix A. As a challenge, we pose the following problem:
\th{Problem.}{}{Find an algorithm that computes the intersection numbers
of the divisor classes $\ll,\d0,\d1,\dots,\d{[g/2]}$ on $\mbar g$.}

\no These numbers are known for $g=2$ [Mu 1] and $g=3$ [Fa 1]. Note that
the problem includes the computation of $\ka1^{3g-3}$.

\th{Proposition}{6.}{Denote by $h_g$ the intersection number
$\ll^{2g-1}\cdot[{\overline {\cal H}}_g]_Q$ , where
${\overline {\cal H}}_g$ is the closure in $\mbar g$ of the
hyperelliptic locus. Then
$$\eqalign{
h_1&=
{1\over96}\quad;\cr
h_g&=
{2\over{2g+1}}\sum_{i=1}^{g-1}i(i+1)(g-i)(g-i+1){{2g-2}\choose{2i-1}}
h_ih_{g-i}\qquad\hbox{for}\qquad g\ge2.
\cr
}
$$}

\Proof This follows from [C-H], Proposition 4.7, which
expresses $\ll$ on ${\overline {\cal H}}_g$ in terms of the classes
of the components of the boundary ${\overline {\cal H}}_g-{\cal H}_g$.
It is easy to see that $\ll^{2g-2}\xi_i=0$ for $0\le i\le[(g-1)/2]$.
Also,
$$
\ll^{2g-2}\d{j}[{\overline {\cal H}}_g]_Q
=(2j+2)(2g-2j+2){{2g-2}\choose{2j-1}}h_jh_{g-j}\quad,
$$
because $\ll=\pi_j^*\ll+\pi_{g-j}^*\ll$ on $\D{j}\cap{\overline {\cal H}}_g$.
Normalizing $h_1$ to ${1\over96}$, which reflects the identity
$\ll={1\over24}p$ on $\mbar{1,1}$ and the fact that an elliptic curve
has four 2-torsion points, we get the formula.
\vskip4pt plus2pt

\pn This gives for instance $h_2={1\over2880}$, $h_3={1\over10080}$
and $h_4={31\over362880}$.
So this already gives the value of $\ll^3$ on $\mbar2$, and the value of
$\ll^6$ on $\mbar3$ follows very easily: we only need that
$[{\cal H}_3]_Q=9\ll$ in $A^1({\cal M}_3)$, because
clearly $\ll^5\d0=\ll^5\d1=0$. We get $\ll^6={1\over90720}$.

\th{Proposition}{7.}{$\ll^9={1\over113400}$ on $\mbar4$ .}

\Proof We need to know the class $[{\overline {\cal H}}_4]$
modulo the kernel in $A^2(\mbar4)$ of multiplication with $\ll^7$.
We computed this class using the test surfaces of [Fa 2]; of the 14
classes at the bottom of p.~432, only $\ka2$, $\ll^2$ and $\d1^2$ are
not in the kernel of $\cdot\ll^7$, and the result is:
$$
[{\overline {\cal H}}_4]\equiv\textstyle3\ka2-15\ll^2+{27\over5}\d1^2
\pmod{\ker(\cdot\ll^7)}.
$$
We also have the relation ([Fa 2], p.~440)
$$
60\ka2-810\ll^2+24\d1^2\equiv0\pmod{\ker(\cdot\ll^7)}.
$$
Thus $[{\overline {\cal H}}_4]\equiv{51\over2}\ll^2+{21\over5}\d1^2$.
We compute
$$\eqalign{\textstyle
\ll^7\d1^2
&=\textstyle{7\choose1}(\ll\cdot[\mbar{1,1}]_Q)
(\ll^6\cdot(-K_{\mbar{3,1}/\mbar3})\cdot[\mbar{3,1}])\cr
&=\textstyle7\cdot{1\over24}\cdot{-4\over90720}\cr
&=\textstyle{-1\over77760}\quad.
}
$$
Therefore
$$
\ll^9=\textstyle{2\over51}(2\cdot{31\over362880}
+{21\over5}\cdot{1\over77760})={1\over113400}\quad.
$$
Also
$$
\ll^7\ka2=\textstyle{169\over1360800}\quad.
$$
The hardest part of this proof is the computation of (three of) the
coefficients of the class $[{\overline {\cal H}}_4]$. We present
the test surfaces we need to compute these coefficients. Write
$$\eqalign{
[{\overline {\cal H}}_4]&=3\ka2-15\ll^2+c\ll\d0+d\ll\d1+e\d0^2+f\d0\d1
+g\d0\d2\cr
&\qquad+h\d1^2+i\d1\d2+j\d2^2+k\d{00}+l\d{01a}+m\gamma_1+n\d{11}\quad.\cr
}
$$
The class $[{\cal H}_4]\in A^2({\cal M}_4)$ was computed by Mumford
([Mu 1], p.~314).
\item{a)} Take test surface $(\alpha)$ from [Fa 2], p.~433: two curves
of genus 2 attached in one point; on both curves the point varies.
We have $[{\overline {\cal H}}_4]_Q=6\cdot6=36$ and $\d2^2=8$. Thus
$j=9$.
\item{b)} Test surface $(\zeta)$: curves of type $\d{12}$, vary the
elliptic tail and the point on the curve of genus 2. We have
$[{\overline {\cal H}}_4]=0$, $\d0\d2=-24$ and $\d1\d2=2$. Thus $i=12g$.
\item{c)} Test surface $(\mu)$: curves of type $\d{02}$, vary the
elliptic curve in a simple pencil with 3 disjoint sections and vary the
point on the curve of genus 2. Then $\d0\d2=-20$ and $\d2^2=4$. To compute
$[{\overline {\cal H}}_4]$ we use a trick. Consider the pencil of curves
of genus 3 which we get by replacing the one-pointed
curve of genus 2 with a fixed one-pointed curve of genus 1. On that pencil
$\ll=1$, $\d0=12-1-1=10$, $\d1=-1$, thus $[{\overline {\cal H}}_3]_Q
=9\ll-\d0-3\d1=2$. So on the test surface we get
$[{\overline {\cal H}}_4]_Q=2\cdot6=12$. Therefore $-20g+36=24$ so
$g={3\over5}$ and $i={36\over5}$.
\item{d)} This test surface is taken from [Fa 3], pp.~72 sq. We take
the universal curve over a pencil of curves of genus 2 as in [A-C], p.~155,
and we attach a fixed one-pointed curve of genus 2. As in [Fa 3] we have
$\ll=3(G-\Sigma)$, $\d0=30(G-\Sigma)$, $\d2=-2G+\Sigma$. Since
$G^2=2$, $G\Sigma=0$ and $\Sigma^2=-2$ we have
$\d0\d2=-60$ and $\d2^2=6$. To compute $\ka2$ we use the same trick
as above: replacing the fixed one-pointed curve of genus 2 by one of genus
1, we get a test surface of curves of genus 3. This will not affect
the computation of $\ka2$; using the formulas of [Fa 1] we find
$\ka2=6$. Also $\d0\Sigma=2\gamma_1$ here, thus $\gamma_1=30$.
Since $[{\overline {\cal H}}_4]=0$, we get
$0=18-60g+6j+30m=30m+36$ so $m=-{6\over5}$.
\item{e)} Test surface $(\lambda)$ from [Fa 2]: curves of type $\d{12}$, vary
both the $j$-invariant of the middle elliptic curve and the (second)
point on it.
We have $\d0\d2=-12$, $\d1\d2=1$, $\d2^2=1$, $\ka2=1$, $\d{01a}=12$ and
$\gamma_1=12$. Since $[{\overline {\cal H}}_4]=0$, we get
$0=3-12g+i+j+12l+12m=12l-{12\over5}$ so $l={1\over5}$.
\item{f)} Test surface $(\kappa)$: curves of type $\d{12}$,
vary a point on the middle elliptic curve and vary the elliptic tail.
Then $\d0\d2=-12$, $\d1\d2=1$, $\d{01a}=-12$, $\d{11}=-1$.
Since $[{\overline {\cal H}}_4]=0$, we find $0=-12g+i-12l-n$ so
$n=-{12\over5}$.
\item{g)} The final test surface we need is $(\gamma)$ from [Fa 2]:
we attach fixed elliptic tails to two varying points on a curve of genus 2.
Then $\d1^2=16$, $\d2^2=-2$, $\ka2=2$, $\d{11}=6$. When the two varying
points are distinct Weierstrass points, we get hyperelliptic curves.
So $[{\overline {\cal H}}_4]_Q=6\cdot5=30$ and we get
$60=6+16h-2j+6n=16h-{132\over5}$ so $h={27\over5}$, as claimed.

\noindent This finishes the proof of Proposition 7.
\vskip4pt plus2pt

\no We can now evaluate the contribution from the constant maps for
$g=4$ (cf.~[BCOV], \S5.13, (5.54)):

\th{Corollary}{8.}{$\la3^3={1\over43545600}$ on $\mbar4$.}

\Proof As explained in [Mu 1], \S5, we have on $\mbar4$ the identity
$$
(1+\la1+\la2+\la3+\la4)(1-\la1+\la2-\la3+\la4)=1\quad.\leqno(*)
$$
One checks that this implies $\la3^3={1\over384}\la1^9$, which finishes
the proof.
\vskip4pt plus2pt

\th{Corollary}{9 (Schottky, Igusa).}{The class of ${\cal M}_4$
in ${\cal A}_4$ equals $8\ll$.}

\Proof
Since $(*)$ holds also on the toroidal compactification
$\widetilde{{\cal A}}_4$, we get
$\la1^{10}=384\la1\la3^3=768\la1\la2\la3\la4$. But it follows from
Hirzebruch's proportionality theorem [Hi 1, 2] that
$$\la1\la2\la3\la4=\prod_{i=1}^4{{|B_{2i}|}\over{4i}}=
{1\over1393459200}\quad,
$$
hence $\ll^{10}={1\over1814400}$ on $\widetilde{{\cal A}}_4$.
Using Theorem 1.5 in [Mu 2] we see that the class of ${\cal M}_4$
in ${\cal A}_4$ is a multiple of $\ll$. Denote by
$t:{\cal M}_4\to{\cal A}_4$ the Torelli morphism and denote by
${\cal J}_4$ its image, the locus of Jacobians. Proposition 7 tells us
that $t^*\ll^9={1\over113400}$. Applying $t_*$ we get
$[{\cal J}_4]\cdot\ll^9={1\over113400}$, hence $[{\cal J}_4]=16\ll$, hence
$[{\cal J}_4]_Q=8\ll$, as claimed.
(The subtlety corresponding to the fact that a general curve of genus
$g\ge3$ has only the trivial automorphism, while its Jacobian has
two automorphisms, appears also in computing $\ll^6$ on $\mbar3$ resp.~on
$\widetilde{{\cal A}}_3$: we saw already that $t^*\ll^6={1\over90720}$;
applying $t_*$ we get $[{\cal J}_3]\cdot\ll^6={1\over90720}$; since
$[{\cal J}_3]=2[\widetilde{{\cal A}}_3]_Q$, we get $\ll^6={1\over181440}$,
which is also what one gets using the proportionality theorem.)
\vskip4pt plus2pt

\ref{BCOV}{
\item{[A-C]} E.~Arbarello and M.~Cornalba, {\it The Picard groups of the
moduli spaces of curves}\/, Topology 26 (1987), 153--171.

\item{[BCOV]} M.~Bershadsky, S.~Cecotti, H.~Ooguri and C.~Vafa,
{\it Kodaira-Spencer Theory of Gravity and Exact Results for Quantum String
Amplitudes}\/,
Commun.~Math.~Phys. 165 (1994), 311--428.\hfill

\item{[C-H]} M.~Cornalba and J.~Harris, {\it Divisor classes associated to
families of stable varieties, with applications to the moduli
space of curves}\/, Ann.~scient.~\'Ec.~Norm.~Sup. 21 (1988), 455--475.\hfill

\item{[Fa 1]} C.~Faber, {\it Chow rings of moduli spaces
of curves I : The Chow ring of
{${\overline {\cal M}}_3$\ }}\/, Ann.~of Math. 132 (1990), 331--419.\hfill

\item{[Fa 2]} C.~Faber, {\it Chow rings of moduli spaces of curves II :
Some results on the Chow ring of {${\overline {\cal M}}_4$\ }}\/,
Ann.~of Math. 132 (1990), 421--449.\hfill

\item{[Fa 3]} C.~Faber, {\it Some results on the codimension-two Chow group
of the moduli space of curves}\/, in:
Algebraic Curves and Projective Geometry (eds E.~Ballico and C.~Cili\-berto),
Lecture Notes in Mathematics 1389, Springer, 66--75.

\item{[Ha]} R.~Hartshorne, {\it Ample Subvarieties of Algebraic Varieties}\/,
Lecture Notes in Mathematics 156, Springer.\hfill

\item{[Hi 1]} F.~Hirzebruch, {\it Automorphe Formen und der Satz von
Riemann-Roch}\/, Symposium Internacional de Topolog\'\i a Algebraica
(M\'exico 1956), 129--144 = Ges.~Abh., Band I, 345--360.\hfill

\item{[Hi 2]} F.~Hirzebruch, {\it Characteristic numbers of homogeneous
domains}\/, Seminars on analytic functions, vol.~II, IAS, Princeton 1957,
92--104 = Ges.~Abh., Band I, 361--366.\hfill

\item{[Ko]} M.~Kontsevich, {\it Intersection
Theory on the Moduli Space of Curves
and the Matrix Airy Function}\/, Commun.~Math.~Phys. 147 (1992), 1--23.\hfill

\item{[Liu]} Qing Liu, {\it Courbes stables de genre 2 et leur sch\'ema de
modules}\/, Math.~Ann. 295 (1993), 201--222.\hfill

\item{[Mu 1]} D.~Mumford, {\it Towards an enumerative geometry of the moduli
space of curves}\/,
in: Arithmetic and Geometry II (eds M.~Artin and J.~Tate),
Progress in Math. 36 (1983), Birkh\"auser, 271--328.

\item{[Mu 2]} D.~Mumford, {\it On the Kodaira Dimension of the Siegel Modular
Variety}\/, in: Algebraic Geometry---Open Problems (eds C.~Ciliberto,
F.~Ghione and F.~Orecchia),
Lecture Notes in Mathematics 997, Springer, 348--375.

\item{[Wi]} E.~Witten, {\it Two dimensional gravity and intersection theory on
moduli space}\/, Surveys in Diff.~Geom. 1 (1991), 243--310.
}

\end